# Exploring iterative and non-iterative Fourier series-based methods of control optimization in application to a discontinuous capsule drive model

Sandra Zarychta, Marek Balcerzak, Jerzy Wojewoda

*Division of Dynamics, Lodz University of Technology, Stefanowskiego 1/15, Lodz, Poland*
sandra.zarychta@dokt.p.lodz.pl, marek.balcerzak.1@p.lodz.pl, jerzy.wojewoda@p.lodz.pl

Abstract

The paper explains iterative and non-iterative approaches to control optimization with use of the Fourier series-based method. Both variants of the presented algorithm are used to numerically approximate optimal control of a discontinuous pendulum capsule drive. Firstly, the general algorithm and its two realizations (iterative and non-iterative) are presented. It is shown that the iterative variant assures non-decreasing quality of solutions in subsequent repetitions of the procedure and the background of such guarantees is explained. A numerical example follows: control of a self-propelled capsule drive is optimized using both approaches. Results are compared and discussed. It is expected that the presented methods can be useful in optimal control estimation for complex systems, particularly discontinuous ones.

*Keywords: iterative, optimal control, Fourier series, open-loop, discontinuous system, self-propelled capsule*

## 1 Introduction

Optimal control of a dynamical system is achieved when it enables a desired action to be performed at the lowest cost. Finding such control is crucial in various scientific and engineering fields [1], [2]. Devices, systems, and processes are expected to operate as planned and produce intended outcomes efficiently. The objective is to achieve goals without wasting time, energy, or other valuable resources, based on some predefined criteria. Consequently, the development of optimal control methods is essential not only from a purely scientific perspective, but also in practical applications. Such control plays a vital role in any area where limited resources need to be utilized efficiently, such as in engineering and economics [1].

In the context of this work, which focuses on optimal control of mechanical systems, it is assumed that the controlled object is described by a set of ordinary differential equations (ODEs) [1], [3] for which determining optimal control usually necessitates the use of numerical methods. The wide range of numerical algorithms for solving Optimal Control Problems (OCPs) can be categorized into three main classes [4], [5]:

- dynamic programming,
- indirect methods based on the calculus of variations,
- direct methods.

The first class comprises methods based on results by Bellman et al.[6], [7], particularly the Hamilton-Jacobi-Bellman (HJB) equation. This approach is practically applicable only for OCPs with a low-dimensional state space, except for the special case of the linear-quadratic regulator (LQR) [1]. The second class includes algorithms grounded in Pontryagin's Minimum Principle (PMP). The OCP is then reduced to solving a boundary value problem [1], [3], typically requiring numerical procedures. Examples of methods in this class include single shooting (parametrizing the entire control function at once), multiple shooting (parametrizing the control function in subsequent time intervals), and collocation methods (where both state and control are parametrized using a predefined space of candidate solutions) [4], [5]. It is important to note that PMP provides necessary but not sufficient conditions for control optimality. The last class encompasses methods based on direct discretization of the OCP, transforming it into a nonlinear programming problem (NLP) to be optimized using a chosen algorithm. Analogously to the previous class, applicable numerical methods include single or multiple shooting, as well as collocation algorithms [4], [5].

Recently, a tool called CASADI has been developed [8], which incorporates methods for nonlinear optimization and algorithmic differentiation in form of a library. It enables effective, algorithmic estimation of derivatives of functions defined as computer program code. Consequently, it offers a comprehensive framework for solving multiple OCPs.

The complexity of OCPs increases when the system to be controlled is defined in terms of a set of discontinuous ODEs. This situation frequently occurs in mechanical systems, where impacts or dry friction play a significant role. In such scenarios, conventional methods for solving optimal control problems (OCPs) are often ineffective [9], [10]. The Pontryagin's Minimum Principle (PMP) requires the vector field to be at least once continuously differentiable with respect to all its arguments [3], making it unsuitable for non-smooth systems. Consequently, standard algorithms based on indirect numerical methods cannot be employed in these problems. Even though the HJB equation does not necessitate differentiation of the vector field, and its smoothness is not formally required, the "optimize the discretization" strategy, fundamental to numerical solutions of HJB equations, leads to excessive errors when the vector field is discontinuous [10].

Several efforts have been made to tackle OCPs involving non-smooth and discontinuous systems. In works [11], [12] OCPs of non-smooth but continuous systems were discussed, while publication [13] addressed OCPs of objects with set-valued vector fields, albeit Lipschitz continuous. Another approach, described in the paper [14], involves computing the Jacobian matrix of the trajectory concerning the initial conditions, feasible as long as neither the initial state nor the final state lies on the discontinuity. According to [14], this method results in conditions analogous to PMP. However, computing this Jacobian matrix is a challenging task. Alternatively, neural networks and evolutionary computation methods have been explored [15]. Lastly, publication [10] presents results involving the approximation of the non-smooth vector field with smooth functions, showing promising outcomes, although it essentially analyzes a qualitatively different system.

A fascinating category of discontinuous control systems is represented by capsule drives, named after their typical capsule-like shapes, which utilize an internal oscillator to generate inertia force. This effect, in combination with dry friction, enables the capsule to move in a desired direction, eliminating the need for external moving parts such as wheels, tracks, or robotic arms. One notable example is the vibro-impact drive [16], where a mass-on-spring oscillator is accompanied by a second spring, allowing the mass to impact it during vibrations. This asymmetrical force results in forward motion of the capsule. Extensive research on this system and its variations has been conducted. The paper [16] presents a comprehensive bifurcation study and discusses optimization of the capsule's motion using simple harmonic control functions. Experimental verification of numerical studies is provided in [17], [18]. Publications [19], [20], [21], [22], [23] analyze various friction models and influence of stochastic conditions. Detailed studies on drive's control, focusing on maximizing progression rate and optimizing energy consumption, can be found in [24], [25], [26], [27], [28]. Multistability control of the vibro-impact drive is detailed in [29], [30], and the design and testing of a small-scale prototype are described in [31]. A comprehensive numerical and experimental research has been conducted with regard to application of the vibro-impact drive in capsule endoscopy, where the system is driven by an external magnetic field [23], [32], [33], [34], [35], [36], [37], [38], [39]. Similar concepts are also applied for monitoring and cleaning tubes or pipelines [40], [41]. Recent achievements in vibro-impact drives development include two-sided motion [20], [42], [43].

Another example is the pendulum capsule drive, where the mass-on-spring oscillator is replaced by a pendulum. Its dynamical properties, analyzed in [44], are more complex due to the pendulum's swinging affecting the contact force between the capsule and the surface, thereby influencing the friction force. Various aspects of this system have been investigated: the impact of viscoelasticity is explored in [45], and friction-induced hysteresis is discussed in [46]. Design and parameter optimization of a pre-defined control function profile are considered in [47], [48]. The controller has been enhanced through the implementation of a neural network-based adaptation mechanism, as described in the paper [49], and energy-related issues related to the control function have been studied in [50]. Construction of the closed-loop, neural controller, with the use of a preoptimized open-loop control, is presented in the work [51]. Application of such machine learning techniques prove particularly effective in systems with complex or unpredictable friction [51], [52], [53] or other significant disturbances [49].

Current research on optimizing capsule drives has generally followed a specific strategy. Initially, a control function's of the "shape" or "profile" is selected and parametrized. Sinusoidal profiles [16], [24], [25], [26] or rectangular waves [31], [35], [36] were commonly used for vibro-impact drives, while object-specific shapes were utilized for pendulum capsule drives [46], [47], [48]. Subsequently, parameters of the chosen profile are optimized. However, it is desired to create methods where the control function's shape is not assumed a priori but is also subject to optimization.

A novel algorithm with such features has been introduced recently [54]. The paper presents a Fourier series-based numerical approach for estimating bounded optimal control. This method is designed to effectively tackle the optimal control problem (OCP) in non-smooth mechanical systems, including those with discontinuities arising from dry friction or impacts. The algorithm employs an innovative technique to parametrize the control function to be optimized, using a finite number of harmonics from the Fourier series. Research indicates that this method significantly enhances the efficiency of the pendulum capsule drive. Similar approach has been applied to optimization of the attitude control of a nanosatellite [55]. Moreover, the aforementioned algorithm has been used as a foundation for closed-loop, neural network based controller [51].

Although the Fourier series-based method has already depicted its undeniable assets, it also has an important drawback: limited scalability. Theoretically, accuracy of the algorithm can be increased by incrementing number of harmonics in the optimized control function. However, this is connected with increasing complexity of the optimization problem under consideration. Taking into account the fact that this approach utilizes evolutionary computation methods such as Differential Evolution [56] or Particle Swarm Optimization [57], which do not guarantee finding the global minimum of the cost functional, it is possible that incrementing number of harmonics in optimized function, intended to improve quality of the solution, may make it worse instead. In the current paper, this problem is solved by introducing a novel, iterative approach. It guarantees that the solution obtained after incrementing the number of harmonics in the optimized control function yields smaller, or at the worst equal, value of the cost functional (which corresponds to non-decreasing control quality).

The paper is organized as follows. Firstly, the Fourier series-based method of control optimization is shortly introduced. Secondly, its modification, which guarantees that quality of solutions is non-decreasing with the number of harmonics, is explained. Further on, the practical example is presented – both approaches are used in optimization of the capsule drive control. Results are compared and the conclusions are drawn.

Authors expect that the presented paper will facilitate research in the field of optimal control, particularly in the field of mechanical systems, in which discontinuities induced by impacts, dry friction and other phenomena are common.

## 2 Description of the method

The purpose of this paper is to develop an incremental, Fourier series-based method of control optimization, in which results of earlier iterations can be used in subsequent trials. Such approach is expected to increase the probability of finding solutions that approach the optimal one. Therefore, this section starts with a precise problem statement, followed by fundamental information on the Fourier series-based method of control optimization [54]. After that, transformation of the algorithm into an incremental method is presented.

### 2.1 Problem statement

Suppose that the system to be controlled is described by an ODE in the following form [1].

$$\dot{x}(t) = f[x(t), u(t), t], \quad x(t_0) = x_0 \quad (1)$$

In the equation (1), $t \in [t_0, t_f] \subset \mathbb{R}$ is the time, $x(t) \in \mathbb{R}^n$ is a state vector at the time $t$, $n$ is the order of the system, $u(t) \in \Omega \subset \mathbb{R}^r$ is a control vector of a dimension $r$ at the time $t$, $f: \mathbb{R}^n \times \mathbb{R}^r \times \mathbb{R} \to \mathbb{R}^n$ is a vector field, $\Omega \subset \mathbb{R}^r$ is a bounded set of admissible controls, $x_0 \in \mathbb{R}^n$ is a vector of initial conditions, $t_0$ and $t_f$ are the initial time and the final time of the control process, respectively.

Any piecewise smooth function $u: [t_0, t_f] \to \Omega$ is called an admissible control function, or shortly, an admissible control. It is assumed that for any admissible control $u$, there exists a unique solution (trajectory) $x: [t_0, t_f] \to \mathbb{R}^n$, $x(t_0) = x_0$ of the system (1).

For any trajectory $x: [t_0, t_f] \to \mathbb{R}^n$ of the equation (1), accompanied by the corresponding admissible control $u: [t_0, t_f] \to \mathbb{R}^r$, the cost functional is defined as follows [1], [2].

$$J = h[x(t_f), t_f] + \int_{t_0}^{t_f} g[x(t), u(t), t] dt \quad (2)$$

The scalar functions $h: \mathbb{R}^n \times \mathbb{R} \to \mathbb{R}$ and $g: \mathbb{R}^n \times \mathbb{R}^r \times \mathbb{R} \to \mathbb{R}$ should be selected in such a manner that desirable trajectories and controls of the system (1) yield smaller values of the cost functional (2) than the unwanted ones.

The main goal of this paper is to describe an effective method that numerically solves the optimal control problem, stated as follows: given the system $f$ (1), its initial conditions $x_0$ and the performance measure $J$ (2), find an optimal control $u^*: [t_0, t_f] \to \Omega$, for which $J$ attains a global minimum [1], [2].

In this paper, the set $\Omega$ in the form of a hyperrectangle is analyzed, i.e. for each component $u_i, i \in \{1, ..., r\}$ of the control function $u(t)$, there exist real constants $m_i, M_i$ such that $\forall_{t \in [t_0, t_1]} m_i \leq u_i(t) \leq M_i$. In other words, the set of admissible controls can be described in the form: $\Omega = [m_1, M_1] \times ... \times [m_r, M_r]$.

In the presented approach to control optimization, it is not required that the vector field $f$ is smooth at all points in the state space. However, if the system is non-smooth, it may be necessary to augment the equation (1) with definitions of switching manifolds and corresponding transition laws [58]. Obviously, it is also required that existence and uniqueness of solutions are assured.

## 2.2 Fourier series-based method of control optimization – a short introduction

To make this paper self-contained yet concise, fundamental information on the Fourier series-based algorithm of control optimization is provided. With the use of this method, optimization of a control function can be reduced to a nonlinear programming problem, in which a relatively small number of parameters is optimized with the use of a selected numerical method. For a more extensive description, please refer to the work [54].

Consider a single component $u_j, 1 \leq j \leq r$ of an admissible control $u$. It has been assumed that each component of $u$ is piecewise smooth and bounded. Consequently, $u_j$ is also integrable [59] and can be approximated using a finite number of terms of the Fourier series.

$$\tilde{u}_j(t) = \frac{a_{j0}}{2} + \sum_{k=1}^{K} a_{jk} \cos(k\omega t) + \sum_{k=1}^{K} b_{jk} \sin(k\omega t) \quad (3)$$

In formula (3), $\omega \in \mathbb{R}_+$ is the fundamental frequency, $a_{j0}, a_{jk}, b_{jk} \in \mathbb{R}$, $K \in \mathbb{N}_+$ and $\tilde{u}_j$ is an approximation of $u_j$. If $u_j$ satisfies Dirichlet conditions [59], then in all points, at which $u_j$ is continuous, $\tilde{u}_j$ converges to $u_j$ as $K \to \infty$.

In order to transform the optimal control problem into a nonlinear programming problem, one needs a method to parametrize the expression (3) in such a way, that for the resulting approximation $\tilde{u} = [\tilde{u}_1, ..., \tilde{u}_r]^T$ the admissibility is assured, i.e. $\tilde{u} \in \Omega$. Under this condition, the performance measure $J$ (2) becomes a function of a finite number of parameters, which enables to use nonlinear programming techniques for control optimization. The method of parametrizing the expression (3) is given below.

First of all, the equation (3) can be rewritten in the form (4).

$$\tilde{u}_j(t) = \frac{a_{j0}}{2} + [a_{j1}, b_{j1}, \ldots, a_{jK}, b_{jK}][\cos(\omega t), \sin(\omega t), \ldots, \cos(K\omega t), \sin(K\omega t)]^T =$$
$$= \frac{a_{j0}}{2} + \boldsymbol{H_j}[\cos(\omega t), \sin(\omega t), \ldots, \cos(K\omega t), \sin(K\omega t)]^T \quad (4)$$

The vector $\boldsymbol{H_j} = [a_{j1}, b_{j1}, \ldots, a_{jK}, b_{jK}]$ contains amplitudes of subsequent harmonics of the function $\tilde{u}_j$. If it is not a constant one, the Euclidean norm of $\boldsymbol{H_j}$ is positive: $H_j = |\boldsymbol{H_j}| > 0$ and $\boldsymbol{H_j}$ can be normalized.

$$\overline{\boldsymbol{H}}_j = \frac{\boldsymbol{H_j}}{H_j} = \frac{[a_{j1}, b_{j1}, \ldots, a_{jK}, b_{jK}]}{\sqrt{a_{j1}^2 + b_{j1}^2 + \cdots + a_{jK}^2 + b_{jK}^2}} \quad (5)$$

Then, the formula (4) can be expressed in the final form as:

$$\tilde{u}_j(t) = \frac{a_{j0}}{2} + \overline{\boldsymbol{H}}_j H_j [\cos(\omega t), \sin(\omega t), \ldots, \cos(K\omega t), \sin(K\omega t)]^T \quad (6)$$

It shows that some important properties of the function $u_j$ are connected only with the direction of the vector $\boldsymbol{H_j}$, represented by its normalization $\overline{\boldsymbol{H}}_j$, and the parameters $K, \omega$. For instance, signs of derivatives of $\tilde{u}_j$ of any order, as well as intervals of monotonicity of $\tilde{u}_j$ and its derivatives, depend only on $\overline{\boldsymbol{H}}_j, K$ and $\omega$. They are invariant under any change of $a_{j0} \in \mathbb{R}$ and $H_j \in \mathbb{R}_+$. These properties, which depend on $\overline{\boldsymbol{H}}_j, K, \omega$, but are not affected by $a_{j0}, H_j$, will be jointly referred to as the *shape* of the function $\tilde{u}_j$. Note that neither adding a constant to the expression (6), nor multiplying it by a positive, real number, changes its shape.

The influence of the parameters $a_{j0}, H_j$ on the function $\tilde{u}_j$ is different. Bearing in mind the graph of $\tilde{u}_j(t)$, changing $a_{j0}$ triggers its shift along the ordinate axis. On the other hand, increasing or decreasing the value $H_j$ cause the graph to expand or shrink respectively (in the direction of the ordinate axis). In other words, the parameter $H_j$ is responsible for scaling the graph of $\tilde{u}_j(t)$. Note that the above-mentioned transformations do not affect the shape of the function under consideration. These properties, connected with shifting or scaling, will be shortly referred to as the *span* of $\tilde{u}_j$.

A deeper and more formal discussion on the notions of the shape and the span can be found in the work [54]. Regarding the scope of this paper, it seems sufficient to say that the shape (described only by $\overline{\boldsymbol{H}}_j, K, \omega$) together with the span (influenced also by $a_{j0}, H_j$) provide complete information on the function $\tilde{u}_j$ under consideration. What is more, the shape of the function does not depend on its span. Therefore, the former can be optimized independently from the latter.

Optimization of $\overline{\boldsymbol{H}}_j$ can be performed with the use of spherical coordinates. Firstly, note that $\boldsymbol{H}_j \in \mathbb{R}^{2K}$. Consider a unit hypersphere $S^{2K-1}$ embedded in $\mathbb{R}^{2K}$, i.e. the set containing these points in $\mathbb{R}^{2K}$, whose distance to the origin equals 1. Obviously, the point indicated by the unit vector $\overline{\boldsymbol{H}}_j$ belongs to this hypersphere. Moreover, location of any point on the unit hypersphere $S^{2K-1}$ can be uniquely defined using $2K-1$ spherical coordinates $\boldsymbol{\varphi}_j = [\varphi_{j1}, \varphi_{j2}, \ldots, \varphi_{j(2K-1)}]^T$, such that $\varphi_{j1}, \varphi_{j2}, \ldots, \varphi_{j(2K-2)} \in [0, \pi]$, $\varphi_{j(2K-1)} \in [0, 2\pi)$ and the following relations hold [60].

$$\overline{H}_{j1} = \cos(\varphi_{j1}) \tag{7a}$$

$$\overline{H}_{j2} = \sin(\varphi_{j1})\cos(\varphi_{j2}) \tag{7b}$$

$$\ldots$$

$$\overline{H}_{j(2K-1)} = \sin(\varphi_{j1})\sin(\varphi_{j2})\ldots\sin(\varphi_{j(2K-2)})\cos(\varphi_{j(2K-1)}) \tag{7c}$$

$$\overline{H}_{j(2K)} = \sin(\varphi_{j1})\sin(\varphi_{j2})\ldots\sin(\varphi_{j(2K-2)})\sin(\varphi_{j(2K-1)}) \tag{7d}$$

Consequently, the direction $\overline{\boldsymbol{H}}_j$ of the vector $\boldsymbol{H}_j$ can be unambiguously defined in terms of $2K-1$ angles: $\varphi_{j1}, \varphi_{j2}, \ldots, \varphi_{j(2K-2)} \in [0, \pi]$, $\varphi_{j(2K-1)} \in [0, 2\pi)$. However, the shape of $\tilde{u}_j$ is influenced by two more parameters: $\omega$ and $K$.

The expression (6) can approximate an arbitrary periodic function of a period $T = 2\pi/\omega$ that satisfies Dirichlet conditions. Therefore, by setting $T = t_f - t_0$ and $\omega = 2\pi/T = 2\pi/(t_f - t_0)$ in the equation (6), one can obtain approximation of an arbitrary Dirichlet function over the time interval $[t_0, t_f]$. This suggests that the parameter $\omega$ can be a priori assigned the value $2\pi/(t_f - t_0)$ and does need to be optimized. However, practically its optimization may be beneficial in some cases, for instance, when it is expected that the optimal control is periodic with a period much smaller than $t_f - t_0$.

The parameter $K$ decides on the number of harmonics in the expression (6). Consequently, its increase can improve the accuracy of solution by the cost of optimization time. In the incremental approach, the right value of $K$ appears naturally in the computation process, as it is explained in the next section.

After discussing the parameters $\overline{H}_j, K, \omega$, which govern the shape of the function (6), the values $a_{j0}, H_j$, influencing the span (shifting and scaling), require analysis. Note that, according to the adopted assumptions, an allowable function $\tilde{u}_j$ has to fulfill the condition $\tilde{u}_j(t) \in [m_j, M_j]$ for any $t \in [t_0, t_f]$. First, assume that $\bar{u}_j$ is the function $\tilde{u}_j$ after normalization to the interval $[0,1]$. Precisely, it can be defined by the following relations.

$$\bar{u}_j(t) = \frac{\alpha_j}{2} + \overline{H}_j \beta_j [\cos(\omega t), \sin(\omega t), \ldots, \cos(K\omega t), \sin(K\omega t)]^T \qquad (8a)$$

$$\min_{t \in [t_0, t_f]} \bar{u}_j(t) = 0, \qquad \max_{t \in [t_0, t_f]} \bar{u}_j(t) = 1 \qquad (8b)$$

In the expression (8a), the values $\alpha_j \in \mathbb{R}, \beta_j \in \mathbb{R}_+$ are selected in such a manner that the equations (8b) hold, which makes them unique. Note that the shape of $\tilde{u}_j$ depends on $\overline{H}_j, \omega$ and $K$ only. Therefore, $\bar{u}_j$ and $\tilde{u}_j$ are of the same shape. Now, the function $\tilde{u}_j$ can be expressed in the following form.

$$\tilde{u}_j(t) = m_j + p(1-q)(M_j - m_j) + \bar{u}_j(t)(M_j - m_j)pq, \qquad p,q \in (0,1] \qquad (9)$$

One can easily notice that the maximum value of the expression (9) is $m_j + p(M_j - m_j)$, which is guaranteed to belong to the interval $(m_j, M_j]$ as long as $p \in (0,1]$. Moreover, the minimum of (9) equals $m_j + p(M_j - m_j)(1-q)$, which is also surely in $(m_j, M_j]$ as long as $p, q \in (0,1]$. Consequently, the parameters $p, q \in (0,1]$ precisely describe the shape (i.e., shifting and scaling) of the function $\tilde{u}_j(t)$ in such a way, that any subinterval of $[m_j, M_j]$ can be selected as the range of $\tilde{u}_j$.

Summing up, it has been shown that the shape of $\tilde{u}_j$ can be uniquely expressed in terms of the angles $\varphi_{j1}, \varphi_{j2}, \ldots, \varphi_{j(2K-2)} \in [0, \pi]$, $\varphi_{j(2K-1)} \in [0, 2\pi)$ and the constants $\omega, K$ – see the equations (6) and (7). Then, the range of $\tilde{u}_j$ is uniquely specified by parameters $p, q \in (0,1]$, as presented in the formulas (8), (9). All these values together yield a complete description of $\tilde{u}_j$, which leads to the value of the cost $J$.

$$J = J(\varphi_{j1}, \varphi_{j2}, \ldots, \varphi_{j(2K-2)}, \omega, K, p, q) \qquad (10)$$

Note that, although the quantity $J$ has been defined as a functional (2), in the expression (10) it becomes a function of a finite number of parameters. With the use of the presented form, the value $J$ can be optimized by means of nonlinear programming techniques. Note that the function (10) is not necessarily continuous with respect to all its arguments. Because of that, special method need to be used, such as the Differential Evolution (DE) [56] or the Particle Swarm Optimization (PSO) [57] algorithms. Unfortunately, their computation costs are relatively high and they do not guarantee finding the global minimum. Because of that, it is crucial to arrange the optimization process in a way that supports obtaining best results. The proposed method of accomplishing this goal is presented in the next section.

## 2.3 Fourier series-based method of control optimization – iterative approach

The description presented in the previous section leaves two important, open questions. First of all, it is not clear how to select the value of the parameter $K$. Secondly, it is not known whether results obtained with smaller values of $K$ can be useful in optimization after incrementing that parameter.

In this section, both of these problems are going to be solved. The goal is to create an algorithm, in which the process of control optimization starts with a small value of $K$, which is connected with a fast convergence of DE or PSO methods. Subsequent iterations of the procedure, with the use of higher values of $K$, are going to utilize previous optimal results as the initial conditions. The process will be repeated as long as an improvement in the results is observed.

Suppose that the control has been optimized with the use of the algorithm presented in the previous section for a fixed number of harmonics $K$ in the control function (3). It means that the optimal angles $\boldsymbol{\varphi}_j = [\varphi_{j1}, \varphi_{j2}, \ldots, \varphi_{j(2K-1)}]^T$ and the corresponding parameters $p, q$ have been approximated. Now, the following question arises: if the accuracy of optimal control's approximation is going to be increased by changing the number of harmonics in the expression (3) to $(K+1)$, then what initial conditions of the optimization procedure should be selected in order to start the new optimization with the best result from the previous execution of the optimizer?

The answer is as follows. By inspection of the formula (9), one can notice that the resulting control approximation $\tilde{u}_j$ will remain unchanged after increasing the number of harmonics from $K$ to $(K+1)$ if the values $p, q$ and the function $\bar{u}_j$ (8) are not modified. On the other hand, the expression $\bar{u}_j$ (8) will not be influenced by changing $K$ to $(K+1)$ only if the amplitudes of the last harmonic are equal 0, i.e. $\bar{H}_{j(2K+1)} = \bar{H}_{j(2K+2)} = 0$, and all the remaining values $\bar{H}_{j1}, \ldots, \bar{H}_{j(2K)}$ do not alter. Increasing the number of harmonics by one triggers appearance of two new angles to be optimized: $\varphi_{j(2K)}, \varphi_{j(2K+1)}$. Formulas (7a-7d) show that the new angles influence the values $\bar{H}_{j(2K-1)}, \bar{H}_{j(2K)}$ and $\bar{H}_{j(2K+1)}$. Moreover, one has to take into account the fact that all the optimized angles should belong to the interval $[0, \pi]$, except for the last one, contained in the set $[0, 2\pi)$. Consequently, with $K$ harmonics, $\varphi_{j(2K-1)}$ was the last optimized angle, i.e. $\varphi_{j(2K-1)} \in [0, 2\pi)$, whereas after introducing two new angles, the intervals are as follows: $\varphi_{j1}, \varphi_{j2}, \ldots, \varphi_{j(2K-1)}, \varphi_{j(2K)} \in [0, \pi]$, $\varphi_{j(2K+1)} \in [0, 2\pi)$. Before taking into account the parameters $\varphi_{j(2K)}, \varphi_{j(2K+1)}$, the relations (7c), (7d) held:

$$\bar{H}_{j(2K-1)} = \sin(\varphi_{j1}) \sin(\varphi_{j2}) \ldots \sin(\varphi_{j(2K-2)}) \cos(\varphi_{j(2K-1)})$$
$$\bar{H}_{j(2K)} = \sin(\varphi_{j1}) \sin(\varphi_{j2}) \ldots \sin(\varphi_{j(2K-2)}) \sin(\varphi_{j(2K-1)})$$

and the values $\bar{H}_{j(2K)}, \bar{H}_{j(2K+1)}$ did not exist. Afterwards, the following formulas are valid.

$$\bar{H}_{j(2K-1)} = \sin(\varphi_{j1}) \sin(\varphi_{j2}) \ldots \sin(\varphi_{j(2K-2)}) \cos(\varphi_{j(2K-1)}) \quad (11a)$$
$$\bar{H}_{j(2K)} = \sin(\varphi_{j1}) \sin(\varphi_{j2}) \ldots \sin(\varphi_{j(2K-1)}) \cos(\varphi_{j(2K)}) \quad (11b)$$
$$\bar{H}_{j(2K+1)} = \sin(\varphi_{j1}) \sin(\varphi_{j2}) \ldots \sin(\varphi_{j(2K)}) \cos(\varphi_{j(2K+1)}) \quad (11c)$$
$$\bar{H}_{j(2K+2)} = \sin(\varphi_{j1}) \sin(\varphi_{j2}) \ldots \sin(\varphi_{j(2K)}) \sin(\varphi_{j(2K+1)}) \quad (11d)$$

Comparing the equation (7c) with (11a) may cause a false impression that nothing changes with regard to $\bar{H}_{j(2K-1)}$. However, this is not true: beforehand the angle $\varphi_{j(2K-1)}$ belonged to the interval $[0, 2\pi)$, whereas afterwards it is contained in the set $[0, \pi]$. Therefore, to maintain equality between the expressions (7c) and (11a), one has to keep the value $\varphi_{j(2K-1)}$ unchanged if $\varphi_{j(2K-1)} \leq \pi$, or replace it with $2\pi - \varphi_{j(2K-1)}$ otherwise, which assures that $\cos(\varphi_{j(2K-1)})$ does not change due to restriction of the angle's interval.

Then, the expression (7d) has to be equal (11b). If the angle $\varphi_{j(2K-1)}$ remains unchanged, then $\sin(\varphi_{j(2K-1)})$ is not varied either. In such case, taking $\varphi_{j(2K)} = 0$ is sufficient to assure the equality between (7d) and (11b). Otherwise, if $\varphi_{j(2K-1)}$ is replaced with $2\pi - \varphi_{j(2K-1)}$, then $\sin(\varphi_{j(2K-1)})$ changes sign. Consequently, taking $\varphi_{j(2K)} = \pi$ makes (7d) and (11b) equal after the replacement of $\varphi_{j(2K-1)}$.

Now, regardless of whether $\varphi_{j(2K)} = 0$ or $\varphi_{j(2K)} = \pi$, the equality $\overline{H}_{j(2K+1)} = 0$ holds naturally – see the formula (11c). Finally, in order to have $\overline{H}_{j(2K+2)} = 0$, one can assume $\varphi_{j(2K+1)} = 0$, according to the formula (11d).

In conclusion, if the accuracy of optimal control's approximation is going to be increased by changing the number of harmonics in the expression (3) from $K$ to $(K + 1)$, then to start the new optimization from the previously optimized result, it is sufficient to keep the angles $\varphi_{j1}, \varphi_{j2}, \ldots, \varphi_{j(2K-2)}$ with the corresponding parameters $p, q$ unchanged, transform $\varphi_{j(2K-1)}$ to $2\pi - \varphi_{j(2K-2)}$ if necessary, and use the zero values of the new angles: $\varphi_{j(2K)} = \varphi_{j(2K+1)} = 0$. By doing so, one ensures that the optimization procedure "knows" the previously-optimized result and it is guaranteed that the new optimized function will perform at least as well as the previous one with respect to the assumed cost function $J$ (10).

The observations described above lead naturally to the method of selecting the value $K$ in the optimization procedure. Selection can be performed in the following iterative process, which is also illustrated in the flowchart shown in Fig. 1.

a) Start optimization with a small value $K$, for instance $1, 2$ or $3$.
b) Change the number of harmonics from $K$ to $K + 1$ and repeat the optimization procedure. Use the results from the previous optimization as the initial conditions.
c) If the performance measure $J$ (10) increases, repeat the point b). Otherwise, finish.

By doing so, one is sure that a new iteration of the optimization process utilizes the previous results. On the other hand, it clearly indicates the maximum, reasonable value $K$, such that its further increase makes no sense in the process of optimal control approximation.

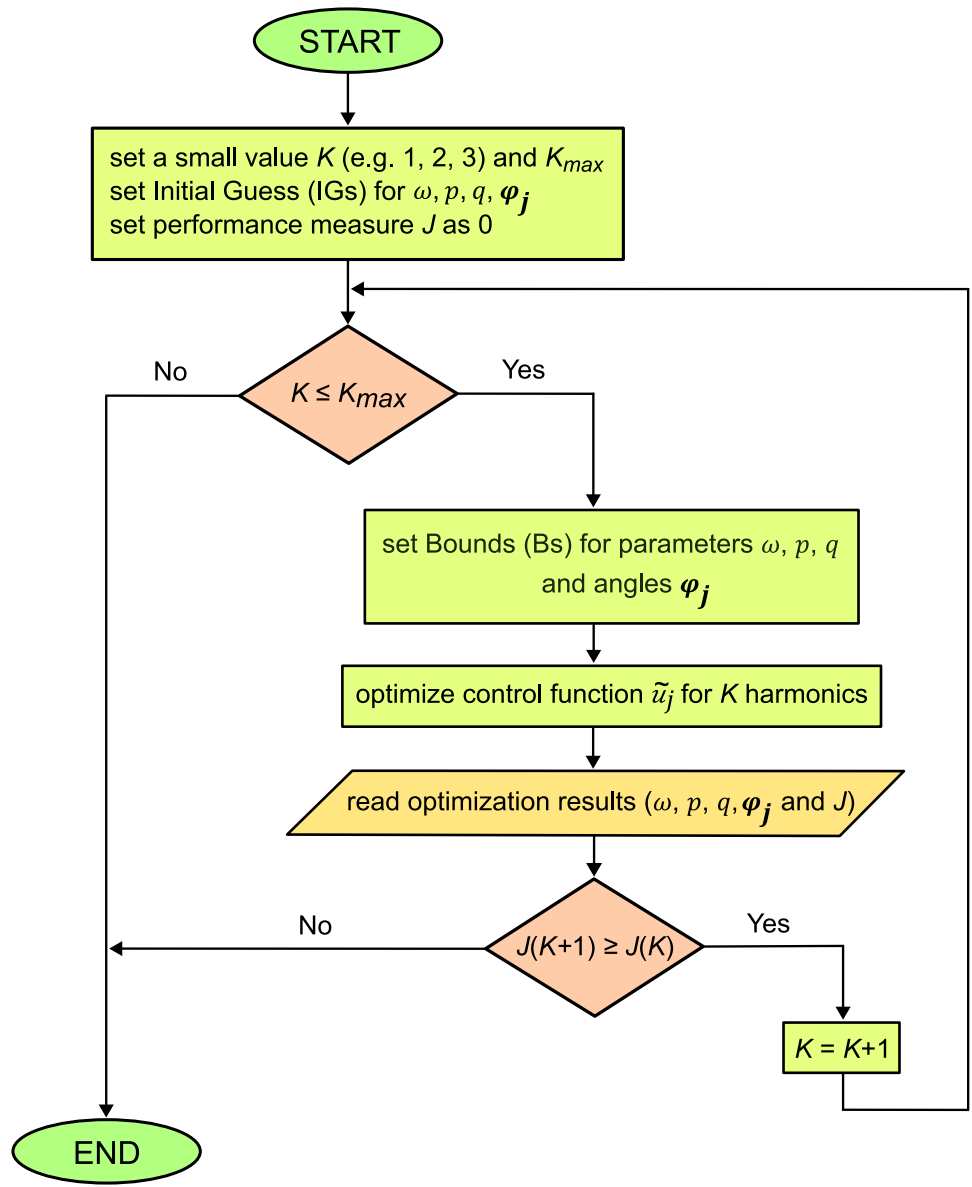

Fig. 1 Flowchart illustrating an iterative approach of the Fourier series-based method of control optimization procedure. $\tilde{u}_j$ – approximation (3) of control function $u_j$; parameters of $\tilde{u}_j$ describing its shape (see exp. 6-7): $\boldsymbol{\varphi}_j$ – angles, $\omega$ - fundamental frequency, $K$ – number of harmonics, $K_{max}$ – maximal number of harmonics in the optimization procedure; $p, q$ (see exp. 8-9) - parameters describing the span of the function $\tilde{u}_j$; $J$ - performance measure (10) of the function $\tilde{u}_j$

# 3 Application of the Fourier series-based methods to a capsule drive model

In the previous section, both approaches of Fourier series-based method of control optimization have been introduced. Firstly, it was shown that a non-iterative algorithm is able to improve the efficiency of discontinuous capsule drive models, albeit with certain restrictions such as limited scalability. Then, the iterative approach was developed to overcome these constraints, ensuring a non-decreasing quality of optimal solutions as the number of harmonics increases.

In this section, both algorithms are applied to a discontinuous capsule drive model and their results are compared. Initially, the equations of the controlled system are formulated. Subsequently, the optimization and simulation details are outlined, followed by a discussion of the resulting findings.

## 3.1 Mathematical model of a discontinuous capsule drive

In this research, control of a pendulum capsule drive is optimized. Such system has been selected in order to enable comparison of results with the non-iterative form of Fourier series based optimization method [54]. Scheme of the device is depicted in Fig. 2. Within this mechanism, the swinging of a mathematical pendulum generates inertia forces that can initiate the motion of the entire capsule, facilitated by the presence of dry friction $F_x$ between the capsule and the surface beneath it. The motion of the pendulum can be induced through the application of an external torque, denoted as $F_\theta(t)$, which needs to be determined.

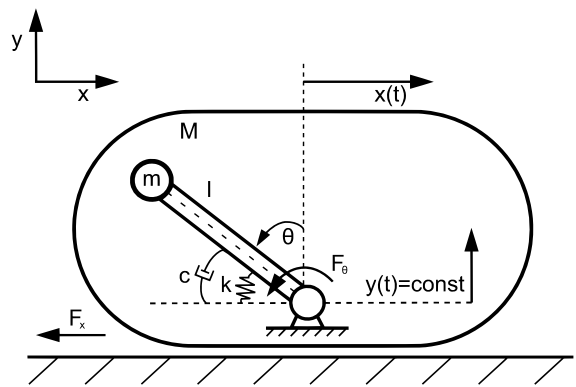

Fig. 2 Scheme of the pendulum capsule drive model. $M$ - mass of the capsule, $m$ – mass of the pendulum, $l$ – length of the pendulum, $\theta$ – pendulum angle, $k$ – spring stiffness, $c$ – damping

coefficient, $F_\theta$ – external torque acting on the pendulum, $F_x$ – friction force, $x(t), y(t)$ – coordinates of the capsule

In the paper [54], a detailed derivation of the motion equations for the pendulum capsule drive have been presented. Therefore, this work provides only a concise overview. The dimensional form of the capsule's motion equations is as follows.

$$ml^2\ddot{\theta}(t) - ml\ddot{x}(t)\cos\theta(t) = mgl\sin\theta(t) - k\theta(t) - c\dot{\theta}(t) + F_\theta(t) \quad (12a)$$

$$(M+m)\ddot{x}(t) - ml\ddot{\theta}(t)\cos\theta(t) + ml\dot{\theta}^2(t)\sin\theta(t) = -F_x(t) \quad (12b)$$

$$R_y(t) = \lambda_y(t) = (M+m)g - ml\ddot{\theta}(t)\sin\theta(t) - ml\dot{\theta}^2(t)\cos\theta(t) \quad (13)$$

Note that $R_y$ represents the vertical, reaction (contact) force between the capsule and the underlying surface. In the system of equations being examined, this force acts as a constraint and is equivalent to the value of the Lagrange multiplier, $\lambda_y$.

The introduced system utilizes the Coulomb friction model, similarly as in the paper [44]. From equation (12a) it can be seen that the horizontal force applied to the capsule due to the pendulum's motion is as follows.

$$R_x(t) = ml\ddot{\theta}(t)\cos\theta(t) - ml\dot{\theta}^2(t)\sin\theta(t) \quad (14)$$

Then, the formulation of the friction model can be established:

$$F_x(t) = \begin{cases} \mu R_y(t) sgn[\dot{x}(t)] \leftrightarrow \dot{x}(t) \neq 0 \\ \mu R_y(t) sgn[R_x(t)] \leftrightarrow \dot{x}(t) = 0 \land |R_x(t)| \geq \mu R_y(t) \\ R_x(t) \leftrightarrow |R_x(t)| < \mu R_y(t) \end{cases} \quad (15)$$

where the following non-dimensional variables and parameters are assumed.

$$\Omega = \sqrt{\frac{g}{l}}, \tau = \Omega t, \gamma = \frac{M}{m}, z = \frac{x}{l}, \rho = \frac{k}{m\Omega^2 l^2}, \nu = \frac{c}{m\Omega l^2},$$
$$f_z = \frac{F_x}{m\Omega^2 l}, u_1 = \frac{F_\theta}{m\Omega^2 l^2}, r_z = \frac{R_x}{m\Omega^2 l}, r_y = \frac{R_y}{m\Omega^2 l} \quad (16)$$

Relations between derivatives regarding the dimensional time $t$ and the dimensionless time $\tau$ are depicted in equation (17).

$$\dot{x} = \frac{dx}{dt} = \frac{dx}{d\tau}\frac{d\tau}{dt} = \Omega\frac{dx}{d\tau} = \Omega x',$$
$$\ddot{x} = \frac{d^2x}{dt^2} = \frac{d}{dt}\left(\frac{dx}{dt}\right) = \frac{d}{d\tau}\left(\Omega\frac{dx}{d\tau}\right)\frac{d\tau}{dt} = \Omega^2\frac{d^2x}{d\tau^2} = \Omega^2 x'' \quad (17)$$

Expressions (16) and (17) enable to represent the equations of motion (12a) and (12b) in a dimensionless, matrix form (18).

$$\begin{bmatrix} 1 & -\cos\theta(\tau) \\ -\cos\theta(\tau) & \gamma+1 \end{bmatrix} \begin{bmatrix} \theta''(\tau) \\ z''(\tau) \end{bmatrix} = \begin{bmatrix} \sin\theta(\tau) - \rho\theta(\tau) - \nu\theta'(\tau) + u_1(\tau) \\ -\theta'^2(\tau)\sin\theta(\tau) - f_z(\tau) \end{bmatrix} \quad (18)$$

In an analogous manner, the definitions of the contact load $R_y$ (13), the resulting horizontal force generated by the pendulum's motion $R_x$ (14), and the friction model $F_x$ (15) undergo a similar transformation.

$$r_y(\tau) = (\gamma+1) - \theta''(\tau)\sin\theta(\tau) - \theta'^2(\tau)\cos\theta(\tau) \quad (19)$$

$$r_z(\tau) = \theta''(\tau)\cos\theta(\tau) - \theta'^2(\tau)\sin\theta(\tau) \quad (20)$$

$$f_z(\tau) = \begin{cases} \mu r_y(\tau) sgn[z'(\tau)] & \leftrightarrow z'(\tau) \neq 0 \\ \mu r_y(\tau) sgn[r_z(\tau)] & \leftrightarrow z'(\tau) = 0 \wedge |r_z(\tau)| \geq \mu r_y(\tau) \\ r_z(\tau) & \leftrightarrow |r_z(\tau)| < \mu r_y(\tau) \end{cases} \quad (21)$$

The equations (18)-(21) establish a comprehensive, dimensionless model of the discontinuous pendulum capsule drive, whose scheme is presented in Fig. 2.

### 3.2 Simulation and optimization details

Here, both the non-iterative and the iterative approaches of Fourier series-based control optimization are applied to the discontinuous capsule drive model (18)-(21). Subsequently, a comparison between them is conducted. In such a manner, it is going to be verified which variant of the method approaches the optimal solution faster, as the number of harmonics $K$ is increased.

In the system under consideration, a single control function is optimized, i.e. the dimensionless torque acting on the pendulum, denoted as $u_j(\tau)$. Obviously, it cannot attain arbitrary values. Therefore, the boundaries for the control function are assumed. It is ensured that at any given time, the condition $u_j(\tau) \in \Omega = [m, M]$ is maintained for constants $m, M \in \mathbb{R}$ such that $m < M$. All simulations were carried out within the same time span from $\tau_0 = 0$ to $\tau_f = 100$, starting from zero initial conditions: $\theta(0) = \theta'(0) = z(0) = z'(0) = 0$ for both methods. This is important due to potential multistability in the controlled system [29]. The assumed values of system's parameters, as well as the constraints for the control function are depicted in Table 1.

Table 1 Summary of the assumed values of the system's parameters (18)-(21) and constraints imposed on the control function $\tilde{u}_j$ (9) for the optimization procedure and simulations

| System's parameters | Value |
|---|---|
| $\mu$ | 0.3 |
| $\rho$ | 2.5 |
| $\nu$ | 1.0 |
| $\gamma$ | 10.0 |
| **Constraints for the control function $\tilde{u}_j$** | |
| lower bound $m$ | -4 |
| upper bound $M$ | 4 |

In the following optimization, applicable for both Fourier series-based algorithms, the objective is to identify an allowable control $u_j: [\tau_0, \tau_f] \to \Omega = [m, M]$ that maximizes the distance covered by the capsule in the dimensionless time interval $[\tau_0, \tau_f]$. Then, the performance measure (10) aimed for minimization can be precisely articulated as follows.

$$J = -|z(\tau_f) - z(\tau_0)| \qquad (22)$$

The system's symmetry leads to the appearance of an absolute value in equation (22). This implies that if a control $u_j(\tau)$ produces distance $z(\tau_f)$, then the control $-u_j(\tau)$ generates $-z(\tau_f)$. Thus, the primary concern is discovering a distance with a large absolute value, while adjusting the motion's direction is feasible by altering the control function's sign, making it a secondary issue.

Following the guidelines presented in subsections **2.2** and **2.3**, the approximation of the control function $\tilde{u}_j$ is uniquely determined by the parameters $\boldsymbol{\varphi}_j, p, q, \omega, K$ – see equations (6)-(9). Each single optimization maintains a fixed number of harmonics, denoted as $K$. Consequently, the cost expression (10) becomes a function $J(\boldsymbol{\varphi}_j, p, q, \omega)$, associating a specific set of parameters $\boldsymbol{\varphi}_j$, $p, q, \omega$ (for a fixed $K$) with the negative absolute distance $-|z(\tau_f) - z(\tau_0)|$ covered by the discontinuous capsule drive model within the assumed time interval. These function values are numerically assessed using the transformation formulas for spherical coordinates (7a)-(7d) to compute the control (6) of a correct shape, yet with a span to be adjusted. Once the function's shape is defined by the spherical coordinates $\boldsymbol{\varphi_j}$, along with parameters $\omega$ and $K$, its range is adapted according to relations (8a)-(8b), with use of formula (9). Subsequently, the obtained control is applied in simulation of the capsule system (18)-(21) within the time interval from $\tau_0 = 0$ to $\tau_f = 100$, resulting in the computation of the final distance $z(\tau_f)$. The parameter ranges and the conditions for optimizing the control function $\tilde{u}_j$ are presented in Table 2. Note that within the iterative approach, as the optimization procedure progresses to subsequent iterations by changing the number of harmonics from $K$ to $K + 1$, the previously obtained results need to satisfy the conditions outlined in expressions (11a)-(11d). This step is crucial to utilize the prior results as the new initial conditions. Throughout this research, 5 trials were performed for each optimization.

Table 2 Parameter ranges for the optimization procedure of control function $\tilde{u}_j$ (9)

| Parameter | Lower bound | Upper bound | Method |
|---|---|---|---|
| $\omega$ | $\dfrac{2\pi}{t_f - t_0}$ | 10.0 | Both (non-iterative and iterative) |
| $p$ | 0.0 | 1.0 | |
| $q$ | 0.0 | 1.0 | |
| $\varphi_{j1}, \varphi_{j2}, \ldots, \varphi_{j(2K-2)}$ | 0.0 | $\pi$ | Both, but in iterative **only before** changing the number of harmonics from $K$ to $K + 1$ – see (7a)-(7d) |
| $\varphi_{j(2K-1)}$ | 0.0 | $2\pi$ | |
| $\varphi_{j1}, \varphi_{j2}, \ldots, \varphi_{j(2K-1)}, \varphi_{j(2K)}$ | 0.0 | $\pi$ | Iterative – after changing the number of harmonics from $K$ to $K + 1$[1] – see (7c)-(7d) and (11a)-11(d) |
| $\varphi_{j(2K+1)}$ | 0.0 | $2\pi$ | |

---

[1] Please note that when initiating a new optimization based on prior results, there might be a need to transform $\varphi_{j(2K-1)}$ to $2\pi - \varphi_{j(2K-2)}$ and to establish new angles $\varphi_{j(2K)} = \varphi_{j(2K+1)} = 0$. For further details please refer to the subchapter **2.3**

In the numerical investigation, the transformation of parameters $\boldsymbol{\varphi}_j, p, q, \omega, K$ into the approximated control function $\tilde{u}_j(\tau)$ (9) is accomplished through a Python script. Subsequently, the Fourier parameters $a_{j0}, a_{jk}, b_{jk}$, where $1 \leq k \leq K$, are passed to a dynamically linked library (DLL) responsible for simulating the system (18)-(21). This library was developed in C++ with the use of the Boost::Odeint library for ODE integration and compiled with the g++ compiler. The simulation employs the RK45 [61] variable step integration method, utilizing an absolute tolerance of $10^{-9}$ and the relative tolerance $10^{-12}$. The identification of stick/slip transitions utilizes the bisection method [62]. Lastly, optimization of the function $J(\boldsymbol{\varphi}_j, p, q, \omega)$, equivalent to optimal control approximation, is carried out using the Differential Evolution method [56] implemented in the SciPy package for Python. All scripts related to the optimization procedure, simulations in Python, and the developed C++ library are available in the research data linked to this paper [63].

**3.3 Results and discussion**

The average maximal distance $z(\tau_f)$ covered by the discontinuous pendulum capsule drive model within the non-dimensional time span $0 \leq \tau \leq 100$ in 5 trials is illustrated in Fig. 3. The variations in the maximum distance depend on the number of harmonics $K$ utilized in the optimization process. Fig. 3 depicts that while the non-iterative approach optimizing control function $u_j$ cannot ensure the distance improvement when transitioning from $K$ to $K + 1$ number of harmonics, the iterative approach consistently enhances solution quality.

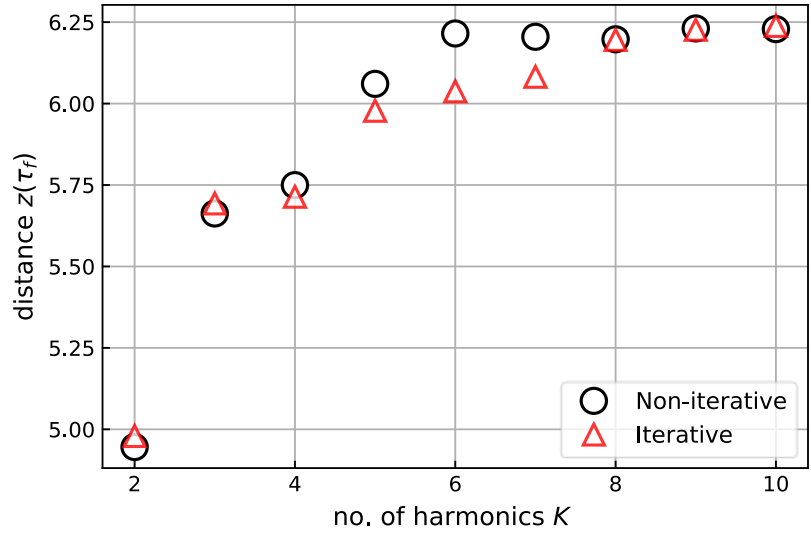

Fig. 3 Maximal distance covered by the discontinuous pendulum capsule drive $z(\tau_f)$ in relation to the number of harmonics $K$ reveals the increasing quality of solution for the iterative approach, whereas it cannot be assured for the non-iterative one (see $K = 7, 8, 10$). Black dots – non-iterative approach; red dots - iterative approach.

To provide further insight into the functioning of both the non-iterative and the iterative approaches in application to the capsule drive, the obtained distances shown in Fig. 3 were used to calculate the relative changes for each method separately, using formula (23):

$$\delta = \frac{z_{K+1} - z_K}{z_K} \cdot 100\% \qquad (23)$$

where:

$\delta$ – the relative change of the distance covered due to the transition from
 $K$ to $K + 1$ harmonics (%),

$z_K$ – the distance covered by the capsule system for $K$ harmonics,

$z_{K+1}$ – the distance covered by the capsule system for $K + 1$ harmonics.

In Table 3, a summary of the achieved results is contained.

Table 3 Summary of the results for non-iterative and iterative approaches of the Fourier series-based method of control optimization: maximal distance covered by the capsule $z(\tau_f)$ along with the standard deviation (SD) depending on the number of harmonics $K$, relative change of distance.

| $K$ | Non-iterative | | Iterative | |
|---|---|---|---|---|
| | Distance $z(\tau_f)$ (±SD) | Relative change [%] | Distance $z(\tau_f)$ (±SD) | Relative change [%] |
| 2 | 4.946 ± 0.029 | - | 4.979 ± 0.019 | - |
| 3 | 5.663 ± 0.026 | ↑ 14.50 | 5.693 ± 0.032 | ↑ 14.34 |
| 4 | 5.750 ± 0.032 | ↑ 1.54 | 5.714 ± 0.044 | ↑ 0.36 |
| 5 | 6.061 ± 0.028 | ↑ 5.41 | 5.978 ± 0.158 | ↑ 4.62 |
| 6 | 6.215 ± 0.045 | ↑ 2.55 | 6.038 ± 0.198 | ↑ 1.00 |
| 7 | 6.205 ± 0.072 | ↓ 0.16 | 6.082 ± 0.161 | ↑ 0.74 |
| 8 | 6.197 ± 0.086 | ↓ 0.12 | 6.195 ± 0.128 | ↑ 1.86 |
| 9 | 6.231 ± 0.090 | ↑ 0.54 | 6.226 ± 0.091 | ↑ 0.50 |
| 10 | 6.229 ± 0.039 | ↓ 0.03 | 6.238 ± 0.082 | ↑ 0.19 |

The results presented in Table 3 offers valuable insights into the functionality of Fourier series-based algorithms. Firstly, the outcomes from the non-iterative approach do not provide the continuous improvement in the quality of solutions (distance) when transitioning from $K$ to $K + 1$. This is notably evident when examining the relative change values. For the initial number of harmonics (up to $K = 6$), all trials consistently show an increase in distance value during the transition from $K$ to $K + 1$. When transitioning from $K = 2$ to $K = 3$, the relative change reaches even 14.50%, indicating a significant increase in the distance covered by the system. Crossing $K = 7$, a pattern of alternating ascent and descent appears for subsequent number of harmonics $K$ in the case of the non-iterative method. The optimization procedure ends at $K = 10$, resulting in a decrease in the distance compared to the preceding harmonic. The highest result achieved using this method is at $K = 9$, attaining the value of 6.231. The distribution of values achieved across 5 subsequential trials is presented as a heatmap as is shown in Fig. 4a). To emphasize the changes in quality of solutions, the following colors are used white for –small improvements, green for –significant ones, and yellow decreases.

Consider the trial no. 3 (as an example): the transition from $K = 2$ to $K = 3$ results in an 14.44% improvement in the distance covered by the system. The next three transitions maintain this positive tendency, but at lower levels (1.45%, 4.94%, 3.32%). When approaching $K = 7$, an unexpected decrease of -0.29% is observed. This decay persists for the subsequent number of harmonics, $K = 8$, resulting in an even higher, -3.25% decrease. After crossing $K = 9$, there is a resumption of improvement in the distance (0.55%, 2.44%) comparing to the preceding values $K$.

The outcomes of the non-iterative approach may lead to the question: why increasing the number of harmonics $K$ does not result in a more precise approximation of the optimal control function $\tilde{u}_j$, as it could be anticipated? As discussed in the introduction part and subsection **2.2**, the developed method employs the Differential Evolution (DE) optimization technique, which does not ensure finding the global minimum of the cost function. Additionally, each optimization process starts from random initial conditions. Consequently, the increase in the number of harmonics $K$ does not necessarily yield better quality of solutions (in this case, the distance). In response to this issue, the iterative approach to control optimization is proposed.

Unlike the non-iterative approach, the iterative method demonstrates that changing the number of harmonics from $K$ to $K + 1$ results in a constantly increasing distance covered by the discontinuous pendulum capsule drive during the optimization procedure (see Fig. 3). Examining the values of relative change from Table 3, a notable improvement occurs between $K = 2$ and $K = 3$ (similarly to the non-iterative approach), reaching 14.34%. Subsequent optimizations show enhancements in solution quality ranging from 0.19% to 4.62%. The optimization procedure ends at $K = 10$ with a 0.19% increase comparing to the preceding value $K$, resulting in the highest observed distance value of 6.238. More detailed outcomes from 5 subsequent trials supporting this analysis are presented in Fig. 4b). The changes in the quality of solutions highlighted by the same colors used in the non-iterative approach (white, green, yellow), do not indicate a decrease in solution quality. This is confirmed by the absence of the yellow color in the data.

It can be noticed that, after crossing $K = 3$, in some trials the obtained distances match those of the previous number of harmonics. For instance, in the trial no. 4, transition from $K = 3$ to $K = 4$ yield no distance improvement, the relative change is equal to 0. This pattern continues until $K = 7$, where an observable distance increase 2.31% is noted, followed by a value of 9.24%. At $K = 9$, only a slight improvement in distance (0.03%) from the previous to the new number of harmonics is observed. Finally, in the last optimization step, there is no observable distance increase.

What is worth noting is that as the number of harmonics $K$ increases, especially when reaching higher values, such as $K = \{9,10\}$, significant distance improvement is observed in only one instance. This observation leads to a conjecture that probability of improving the solution in subsequent steps decreases as the number of harmonics is incremented. Consequently, if no improvement is observed after increasing the value $K$ and the control quality seems satisfactory, it may be reasonable to terminate the procedure. Nevertheless, although improvement of the solution is less probable for large $K$, it is usually possible[2], and the iterative approach guaranties that as $K$ increases, the result is at least not deteriorated.

This section of results analysis distinctly addresses the inquiries outlined in subsection **3.2**. The expectation during the optimization process was that increasing the number of harmonics $K$ would yield a more precise approximation of the optimal control function $\tilde{u}_j$. However, the outcomes for the non-iterative approach indicate that changing the number of harmonics from $K$ to $K + 1$ in pendulum's forcing function $u_j$ may not lead to a smaller value of the cost functional in subsequent iterations within the specified dimensionless time interval $\tau$. Meanwhile, in the iterative approach, the results obtained for $K + 1$ harmonics at least match those for $K$. The non-decreasing quality of solutions ensured by the incremental, Fourier series-based approach confirms that the method works more efficiently compared to the non-iterative algorithm.

---

[2] It is very unlikely that the optimal control is exactly described by the formula (3) with a finite number of harmonics. In any other case, increasing the value $K$ can lead to improvement of the approximate solution.

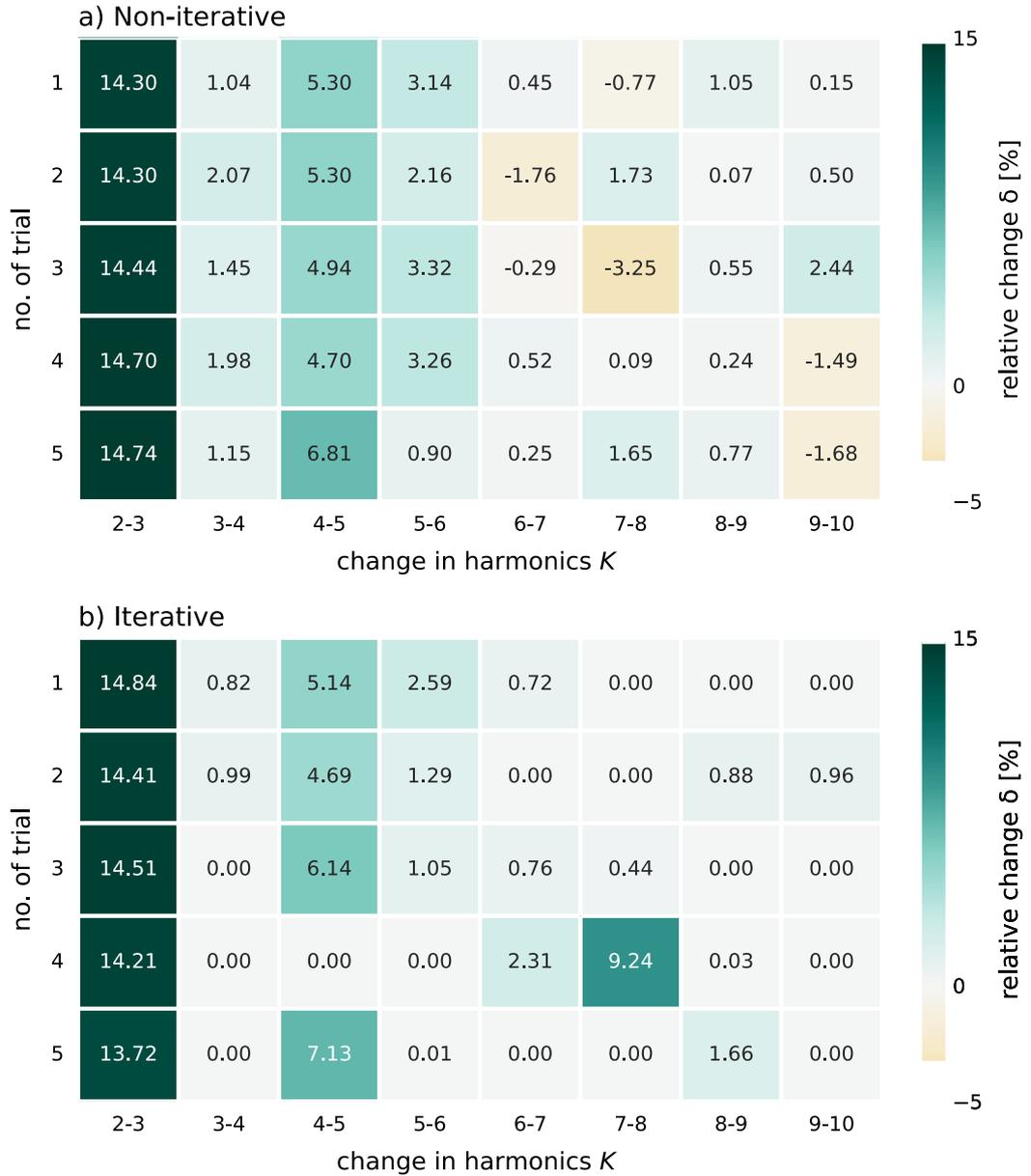

Fig. 4 Heatmaps providing details about relative changes across 5 subsequent trials during the control optimization procedure concerning varying number of harmonics $K$: (a) non-iterative approach – alternating increases and decreases in solution quality (see when $K > 7$); (b) iterative approach – consistently non-decreasing quality of solutions ($\delta \geq 0$).

The highest distances obtained for both, non-iterative and iterative approaches are presented in Fig. 5. The figure visualizes how different numbers of harmonics in the optimized control functions affect the temporal changes in the capsule's position $z(\tau)$. The motion patterns for both methods share similarities – neither follows a monotonic characteristic. This observation suggests that achieving optimal control of the capsule drive involves slight backward movement at certain moments of time to facilitate subsequent accelerated forward motion.

Approximated control functions for the non-iterative method ($K = 9$) and the iterative one ($K = 10$) are depicted in Fig. 6.

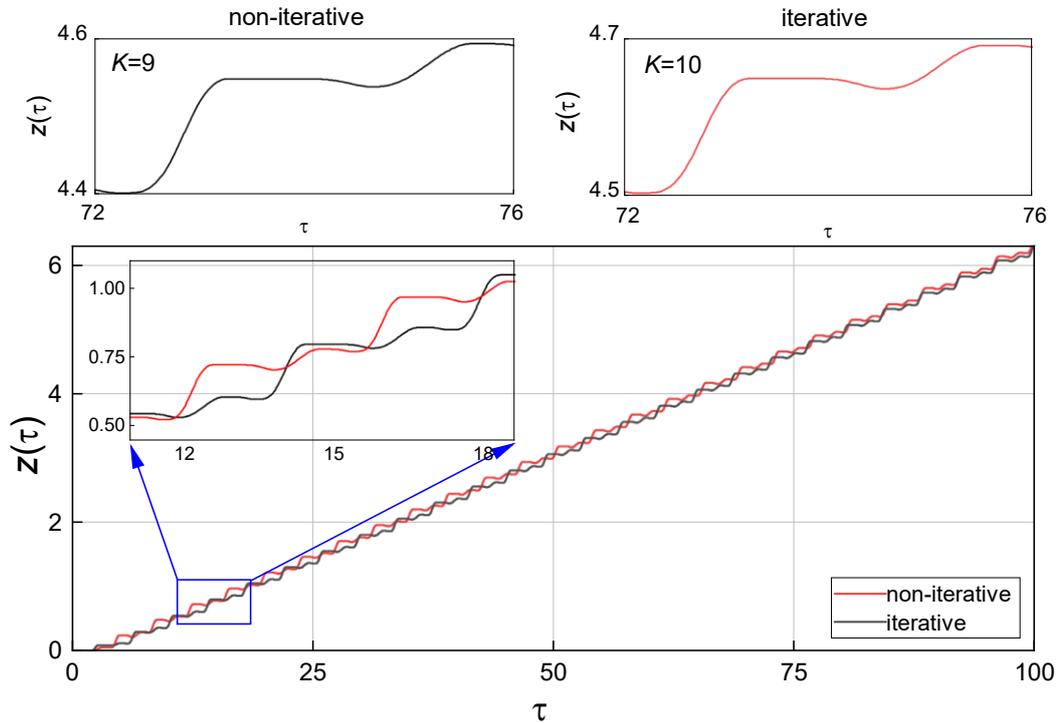

Fig. 5 Position of the capsule $z(\tau)$ vs. dimensionless time $\tau$ for non-iterative and iterative approaches – top results. Changes in position reveal non-monotonic characters of motion, where the stick-slip phenomenon is clearly visible: the flat pieces of graphs (stick) alternate with inclined fragments (slip) in a repetitive manner. Black line – non-iterative approach; red line - iterative approach.

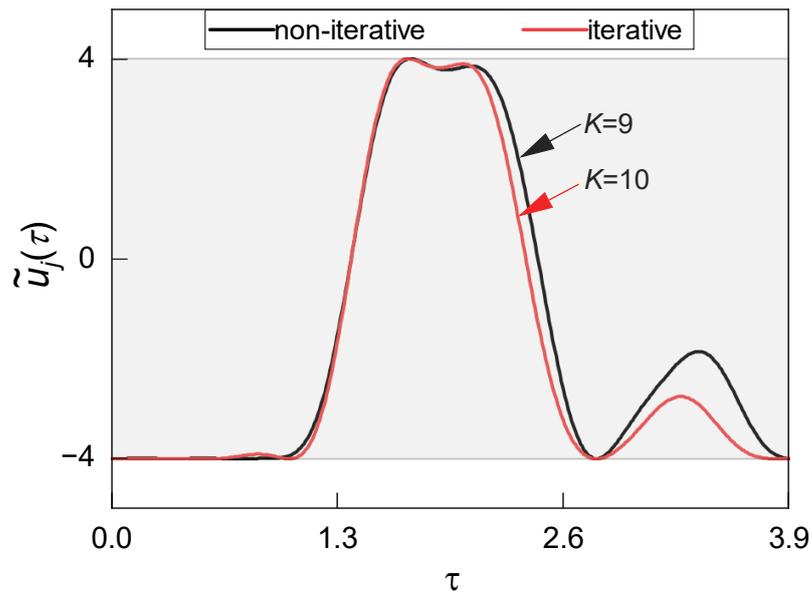

Fig. 6 Optimal control approximations $\tilde{u}_j$ for non-iterative and iterative methods. Black line – non-iterative approach; red line - iterative approach.

# 4 Conclusions

The iterative and non-iterative variants of Fourier series-based method of control optimization in application to a discontinuous system have been explored. Section 2 outlines the problem, introduces the non-iterative approach of the algorithm defining its limitations to solve by the iterative approach, and provides extensive explanation of the incremental optimization algorithm. In Section 3, both variants were applied to optimize the control of a complex discontinuous system: - the pendulum capsule drive. Its motion is initiated by oscillations of an internal mathematical pendulum in the presence of dry friction. Dimensionless parameters are introduced, leading to non-dimensional equations, followed by an explanation of the simulation procedure and presentation of its results.

The described methods offer a straightforward and convenient means of estimating open-loop optimal control while addressing control constraints. These algorithms demonstrate high flexibility, requiring minimal information about the controlled object. The primary prerequisites involve specifying admissible controls within constant ranges and the ability to assess the performance measure $J$ for any admissible control, ensuring its uniqueness. Notably, they can accommodate non-smooth or discontinuous systems, provided their solutions exist and are unique for any admissible control. Additionally, these methods simplify function optimization by transforming it into the more manageable task of parameter optimization (nonlinear programming). This transformation enables solutions through various global optimization algorithms, including the widely-used Differential Evolution procedure.

It has to be noted that although both: the iterative and non-iterative approaches utilize the evolutionary algorithm that cannot guarantee the discovery of the optimal control function or ensure the convergence of the optimization process (they are likely to yield solutions that closely approximate the optimal one), the iterative approach seems overcome this limitation (at least partly). Based on the analyzed numerical example, although the optimality of the derived control is not assured, the solution quality remains non-decreasing as the number of harmonics $K$ changes during the optimization procedure (relative change $\delta \geq 0$). Moreover, as the cost function $J$ does not decrease when compared to the previous iteration, it suggests a reasonable number of harmonics $K$, beyond which further optimization seems futile. The iterative algorithm, in certain scenarios, exceeds the distance achieved in the prior harmonics by over 9%. It is worth noting that the initial transition from $K = 2$ to $K = 3$ yields comparable results for both methods (~14%) and is not explicitly mentioned here. Meanwhile, in the non-iterative approach the absence of improvement becomes evident when $K > 7$. The observations highlight the superiority of the iterative method's efficiency over the non-iterative one in optimizing capsule drives.

While absolute assurance of similar performance across all systems cannot be guaranteed, the attained results so far show considerable promise. The incremental algorithm's potential applications span a wide range, and favorable outcomes are expected across diverse control systems. The authors anticipate that this method can be applied not only to various types of capsule drives, but also enable approximate investigations into optimal control within previously challenging or even deemed impossible areas.

**Acknowledgements**

This study has been supported by the National Science Centre, Poland, PRELUDIUM Programme (Project No. 2020/37/N/ST8/03448).

Work of Marek Balcerzak has been funded under the 'FU2N - Fund for Enhancing the Skills of Young Scientists' program supporting scientific excellence at the Lodz University of Technology – grant No. 55/2021.

This paper has been completed while the first author was the Doctoral Candidate in the Interdisciplinary Doctoral School at the Lodz University of Technology, Poland.

**Data Availability Statement**

The datasets generated and analyzed during the current study are available from the corresponding author on reasonable request. All the scripts created within this study are available in the Mendeley Data repository linked to this paper [63].

**Declarations**

**Conflict of Interest**

The authors declare that they have no conflict of interest concerning the publication of this manuscript.